\documentclass[12pt, reqno]{amsart}
\usepackage{amsmath, amsthm, amscd, amsfonts, amssymb, graphicx, color}
\usepackage[bookmarksnumbered, colorlinks, plainpages]{hyperref}
\hypersetup{colorlinks=true,linkcolor=red, anchorcolor=green, citecolor=cyan, urlcolor=red, filecolor=magenta, pdftoolbar=true}

\newtheorem{theorem}{Theorem}[section]
\newtheorem{lemma}[theorem]{Lemma}
\newtheorem{proposition}[theorem]{Proposition}

\theoremstyle{definition}

\theoremstyle{remark}
\newtheorem{remark}[theorem]{Remark}
\numberwithin{equation}{section}

\begin{document}

\setcounter{page}{1}

\title[Variable Hardy spaces]{Discrete para-product operators
on variable Hardy spaces}

\author[J. Tan]{Jian Tan}

\address{College of Science, Nanjing University of Posts and Telecommunications, Nanjing 210023, China.}
\email{\textcolor[rgb]{0.00,0.00,0.84}{tanjian89@126.com}}





\date{Received: xxxxxx; Revised: yyyyyy; Accepted: zzzzzz.}

\begin{abstract}
Let $p(\cdot):\mathbb R^n\rightarrow(0,\infty)$
be a variable exponent function satisfying the globally
log-H\"older continuous condition.
In this paper, we obtain the boundedness of para-product operators $\pi_b$
on variable Hardy spaces $H^{p(\cdot)}(\mathbb R^n)$, where $b\in BMO(\mathbb R^n)$. As an application,
we show that non-convolution type Calder\'on-Zygmund
operators $T$
are bounded on $H^{p(\cdot)}(\mathbb R^n)$ if and only if $T^\ast1=0$, where
$\frac{n}{n+\epsilon}<\mbox{ess\,inf}_{x\in\mathbb R^n} p\le
\mbox{ess\,sup}_{x\in\mathbb R^n} p\le 1$, $\epsilon$
is the regular exponent of kernel of $T$.
Our approach relies on the discrete version of Calder\'on's reproducing
formula, discrete Littlewood-Paley-Stein theory
and almost orthogonal estimates.
These results still hold for
variable Hardy space on spaces of homogeneous type
by using our methods.
\end{abstract} \maketitle

\section{Introduction and statements of results}

The real-variable theory of Hardy spaces $H^p(\mathbb R^n)$
was initiated by Stein and Weiss
\cite{CW1} and systematically developed by Fefferman and Stein in \cite{FS}.
The Hardy space $H^p(\mathbb R^n)$ with $0<p\leq1$, which is a suitable substitute of the Lebesgue space $L^p(\mathbb R^n)$,
plays an important role in the study of operators
and their application to partial differential equations.

In this paper, we consider the variable Hardy space estimates
for non-convolution type Calder\'on-Zygmund singular
integral operators. Our proofs mainly use the
$H^{p(\cdot)}(\mathbb R^n)-$boundedness
of the discrete para-product operators,
almost orthogonality estimates and
the discrete version of Calder\'on's reproducing
formula. These techniques
and estimates
used in our proofs are
different from the approaches which have been used in the existing
proofs of
many other variable Hardy spaces inequalities.
For example, the boundedness of operators on
variable Hardy spaces $H^{p(\cdot)}(\mathbb R^n)$
was usually established
via using atomic decompositions together with
maximal operator estimate for the action of $T$ on atoms
in recent years.
We refer to the work of Cruz-Uribe and his collaborators
\cite{CMN, CMN1, CW}, Nakai and Sawano \cite{NS,S} and
the author \cite{T2,T3,T4,T5}.
The method to variable Hardy spaces in this paper
can be applied to more general cases even when the maximal function
characterizations are absent.
On the other hand, the para-product operators can be used
to nonlinear analysis. See Bony \cite{B} or Meyer \cite{M}
for more details. There are many different forms of para-products operators,
such as continuous para-product operator \cite{DH1}, frame
para-product operator \cite{DH} and wavelet para-product operator
\cite{HL}. Motivated by these,
we introduce a new discrete para-product operator $\pi_b$
and obtain the $H^{p(\cdot)}(\mathbb R^n)-$boundedness of $\pi_b$,
which also has its own interests.

We first recall some notations and known results on
variable function spaces that will be used in this paper.
See \cite{CF,CW,DHHR,KR,NS,YZN} for more information.
For a measurable subset $E\subset \mathbb{R}^n$, we denote
$p^-(E)= \inf_{x\in E}p(x)$ and $p^+(E)= \sup_{x\in E}p(x).$
Especially, we denote $p^-=p^{-}(\mathbb{R}^n)$ and $p^+=p^{+}(\mathbb{R}^n)$.
Let $p(\cdot)$: $\mathbb{R}^n\rightarrow(0,\infty)$ be a measurable
function with $0<p^-\leq p^+ <\infty$ and $\mathcal{P}^0(\mathbb R^n)$
be the set of all these $p(\cdot)$.
Let $\mathcal{P}(\mathbb R^n)$ be the set of all measurable functions
$p(\cdot):\mathbb{R}^n \rightarrow[1,\infty) $ such that
$1<p^-\leq p^+ <\infty.$
Let $\mathcal{B}(\mathbb R^n)$ be the set of $p(\cdot)\in \mathcal{P}(\mathbb R^n)$ such that the
Hardy-littlewood maximal operator $M$ is bounded on  $L^{p(\cdot)}$.
For any Schwartz functions $\psi$ on $\mathbb R^n$,
we denote that $\psi_t(x)=t^{-n}\psi(t^{-1}x),\;t>0$.
Hereafter, for simplicity we use $\psi_k(x)$
to denote $\psi_{2^{-k}}(x)=2^{kn}\psi(2^{k}x),\; k\in\mathbb N$.

The variable Lebesgue space $L^{p(\cdot)}(\mathbb{R}^n)$
is defined as the set of all
measurable function $f$ for which the quantity
$\int_{\mathbb{R}^n}|\varepsilon f(x)|^{p(x)}dx$ is finite for some
$\varepsilon>0$ and
$$\|f\|_{L^{p(\cdot)}(\mathbb{R}^n)}=\inf{\left\{\lambda>0: \int_{\mathbb{R}^n}\left(\frac{|f(x)|}{\lambda}\right)^{p(x)}dx\leq 1 \right\}}.$$
As a special case of the theory of Nakano and Luxemberg,
we see that $L^{p(\cdot)}(\mathbb{R}^n)$
is a quasi-normed space. Especially, when $p^-\geq1$, $L^{p(\cdot)}(\mathbb{R}^n)$ is a Banach space. In the study of variable exponent function spaces it is common
to assume that the exponent function $p(\cdot)$ satisfies $LH(\mathbb R^n)$
conditions.
We say that $p(\cdot)\in LH(\mathbb R^n)$, if $p(\cdot)$ satisfies

 $$|p(x)-p(y)|\leq \frac{C}{-\log(|x-y|)} ,\quad |x-y| \leq 1/2$$
and
 $$|p(x)-p(y)|\leq \frac{C}{\log|x|+e} ,\quad |y|\geq |x|.$$
It is well known
that $p(\cdot)\in \mathcal{B}(\mathbb R^n)$
if $p(\cdot)\in \mathcal{P}(\mathbb R^n)\cap LH(\mathbb R^n).$
Moreover, examples shows that the above $LH(\mathbb R^n)$ conditions are necessary in certain sense, see Pick and R${\mathring{\mbox{u}}}$\v{z}i\v{c}ka (\cite{PR}) for more details.

Denote by $\mathcal{M}$ the grand maximal operator given by
$\mathcal{M}f(x)= \sup\{|\psi_t\ast f(x)|: t>0,\psi \in \mathcal{F}_N\}$ for any fixed large integer $N$,
where $\mathcal{F}_N=\{\psi \in \mathcal{S}:\int\psi(x)dx=1,\sum_{|\alpha|\leq N}\sup(1+|x|)^N|\partial ^\alpha \psi(x)|\leq 1\}$.
The variable Hardy space ${H}^{p(\cdot)}(\mathbb R^n)$ is the set of all $f\in \mathcal{S}^\prime$ for which the quantity
$$\|f\|_{{H}^{p(\cdot)}(\mathbb R^n)}
=\|\mathcal{M}f\|_{{L}^{p(\cdot)}(\mathbb R^n)}<\infty.$$
Let $\phi \in \mathcal{S}(\mathbb R^n)$ with $\int\phi(x)dx=1$.
Denote that
$$M_\phi(f)(x)=\sup_k|\phi_k\ast f(x)|.$$
If $p(\cdot)\in LH(\mathbb R^n)\cap \mathcal P^0(\mathbb R^n)$,
$$
\|f\|_{H^{p(\cdot)}(\mathbb R^n)}\sim
\|M_\phi(f)\|_{L^{p(\cdot)}(\mathbb R^n)}.
$$

Throughout this paper, $C$ or $c$ denotes a positive constant that may vary at each occurrence
but is independent to the main parameter, and $A\sim B$ means that there are constants
$C_1>0$ and $C_2>0$ independent of the the main parameter such that $C_1B\leq A\leq C_2B$.
Given a measurable set $S\subset \mathbb{R}^n$, $|S|$ denotes the Lebesgue measure and $\chi_S$
means the characteristic function.
For a cube $Q$, let $Q^\ast$ denote with the same center
and $100$ its side length, i.e. $l(Q^\ast)=100l(Q)$.
We define $\mathcal D(\mathbb R^n)$ to be the space of all smooth
functions with compact support and
$\mathcal D_0(\mathbb R^n)$ to be the space of all smooth
functions with compact support and integral zero.
We also use the notations $j\wedge j'=\min\{j,j'\}$.

A local integrable function $b$ is in $BMO(\mathbb R^n)$,
if
$$
\|b\|_{BMO(\mathbb R^n)}=\sup_{Q}\frac{1}{|Q|}\int_Q|f(x)-f_Q|dx<\infty,
$$
where the supremum ranges over all finite cubes $Q\in \mathbb R^n$.
Let $\varphi$ be Schwartz functions with support
on the unit ball satisfying the conditions: for all $\xi\in\mathbb R^n$,
\begin{align*}
\sum_{j\in\mathbb Z}|\widehat\varphi(2^{-j}\xi)|^2=1
\end{align*}
and $\int_{\mathbb R^n}\varphi(x)x^\alpha dx=0$ for
all $0\le|\alpha|\le M$.
For any $b\in BMO(\mathbb R^n)$, the discrete para-product operator $\pi_b$
is defined by
$$
\pi_b(f)(x)=\sum_k\sum_Q|Q|\varphi_k(x-x_Q)\tilde\varphi_k\ast
b(x_Q)\phi_k\ast f(x_Q)
$$
and
$$
\pi^\ast_b(f)(x)=\sum_k\sum_Q|Q|\varphi_k\ast f(x_Q)\tilde\varphi_k\ast
b(x_Q)\phi_k(x_Q-x),
$$
where $x_Q$ is any point in $Q$ and the summation of $Q$ is taken over all dyadic cubes $Q$ with side length $2^{-k-N}$ in $\mathbb{R}^n$ for each $k\in\mathbb Z$ and
a fixed large integer $N$.

First we obtain the $H^{p(\cdot)}(\mathbb R^n)$ boundedness of $\pi_b$.

\begin{theorem}\label{s1th1}\quad
Suppose that $b\in BMO(\mathbb R^n)$ and
$p(\cdot)\in LH(\mathbb R^n)\cap \mathcal P^0(\mathbb R^n)$.
Then $\pi_b$ is bounded on $H^{p(\cdot)}(\mathbb R^n)$
and bounded from $H^{p(\cdot)}(\mathbb R^n)$ to $L^{p(\cdot)}(\mathbb R^n)$.
\end{theorem}

Let $\epsilon$ be the regularity exponent with $0<\epsilon\le 1$.
An operator $T$ is a non-convolution Calder\'on-Zygmund
operators, denoted by $T\in CZO(\epsilon)$,
if it is bounded operator on $L^2(\mathbb R^n)$,
and if for all $f, g\in \mathcal D(\mathbb R^n)$
with $\mbox{supp}(f)\cap\mbox{supp}(g)=\emptyset,$
$$\big<Tf,g\big>=\int_{\mathbb R^n}\int_{\mathbb R^n}K(x,y)f(y)g(x)dydx,
$$
where the distributional kernel coincides with a function
$K$ defined away the diagonal on $\mathbb R^{2n}$
and satisfies the following conditions
\begin{align}
  & |K(x,y)|\le \frac{C}{|x-y|^n}, \label{s1i1}\\
  & |K(x,y)-K(x',y)|\le C\frac{|x-x'|^\epsilon}{|x-y|^{n+\epsilon}},
  \quad \mbox{if}\;|x-x'|\le\frac{1}{2}|x-y|, \label{s1i2}\\
  & |K(x,y)-K(x,y')|\le C\frac{|y-y'|^\epsilon}{|x-y|^{n+\epsilon}},
  \quad \mbox{if}\;|y-y'|\le\frac{1}{2}|x-y|. \label{s1i3}
\end{align}
Obviously, for any $\eta\in \mathcal D_0$,
$$\big<T1,\eta\big>=\int_{\mathbb R^n}\int_{\mathbb R^n}K(x,y)\eta(x)dydx,
$$
and
$$\big<T^\ast1,\eta\big>=\int_{\mathbb R^n}\int_{\mathbb R^n}K(x,y)\eta(y)dydx.
$$

Then we state the following result.
\begin{theorem}\label{s1th2}\quad
Suppose that $p(\cdot)\in LH(\mathbb R^n)$ and $T\in CZO(\epsilon),$ where $0<\epsilon\le 1$.
Then $T$ extends to a continuous operator on $H^{p(\cdot)}(\mathbb R^n)$,$\frac{n}{n+\epsilon}<p^-\le p^+\le1$
if and only if $T^\ast1=0$.
\end{theorem}

\begin{remark}
Theorem \ref{s1th1} and Theorem \ref{s1th2} still hold for
$H^{p(\cdot)}(\mathcal X)$
under some suitable conditions,
where $\mathcal X$ is a space of homogeneous type
introduced by R. Coifman and G. Weiss in \cite{CW1}. We will see
that, in the next sections, our proofs of these theorems
only depends on the atomic decomposition characterization
of variable Hardy spaces, the discrete Calder\'on reproducing formula,
the discrete Littlewood-Paley-Stein theory and
almost orthogonality estimates. The atomic
decomposition theory of $H^{p(\cdot)}(\mathcal X)$
was established in \cite{ZSY}. The formula, the
discrete Littlewood-Paley-Stein theory
and the estimates still hold for the
spaces of homogeneous type. For more details,
see \cite{DH,HMY,HMY1}. Therefore, the theorems still hold
for this setting.
\end{remark}

To study the boundedness of the Calder\'on-Zygmund operators $T$ on
other variable Hardy spaces $H^{p(\cdot)}(\mathbb R^n)$,
we need to assume more additional
conditions on $T$.
We say that an operator $T\in CZO(N,\epsilon)$,
if the distributional kernel coincides with a function
$K$ defined away the diagonal on $\mathbb R^{2n}$
and satisfies the following conditions
\begin{align*}
  & |\partial_x^\alpha K(x,y)|+|\partial_y^\alpha K(x,y)|\le
  \frac{C}{|x-y|^{n+|\alpha|}}, \quad \mbox{if}\;0\le|\alpha|\le N,\\
  & |\partial_x^\alpha K(x,y)-\partial_y^\alpha K(x',y)|\le C\frac{|x-x'|^\epsilon}{|x-y|^{n+N+\epsilon}},
  \quad \mbox{if}\;|\alpha|=N,\; |x-x'|\le\frac{1}{2}|x-y|,\\
  & |\partial_x^\alpha K(x,y)-\partial_y^\alpha K(x,y')|\le C\frac{|y-y'|^\epsilon}{|x-y|^{n+N+\epsilon}},
  \quad \mbox{if}\;|\alpha|=N,\;|y-y'|\le\frac{1}{2}|x-y|,
\end{align*}
where $0<\epsilon\le 1$ and $N>0$ is an integer.

Finally, we obtain the following theorem.
\begin{theorem}\label{s1th3}\quad
Suppose that $p(\cdot)\in LH(\mathbb R^n)$ and $T\in CZO(N,\epsilon)\cap CZO,$ where $0<\epsilon\le 1$.
Then $T$ extends to a continuous operator on $H^{p(\cdot)}(\mathbb R^n)$,$\frac{n}{n+N+\epsilon}<p^-\le p^+\le1$
if and only if $T^\ast x^\beta=0$ with $|\beta|\le N$.
\end{theorem}

\begin{remark}
After we were completing this paper,
we learned that Cruz-Uribe, Moen
and Nguyen also obtained the boundedness of
non-convolution Calder\'on-Zygmund operators
on variable Hardy spaces via the different method, which
mainly contains finite atomic decompositions
and the extrapolation method (See \cite[Theorem 1.11]{CMN1}).
\end{remark}

 \section{Proof of Theorem \ref{s1th1}}
In this section, we will establish the $H^{p(\cdot)}(\mathbb R^n)$
boundedness of $\pi_b$. To do this, we need the atomic decomposition theory for $H^{p(\cdot)}(\mathbb R^n)$.
We recall the atoms of variable Hardy spaces ${H}^{p(\cdot)}(\mathbb R^n)$
in \cite{NS}.
Let $p(\cdot): \mathbb{R}^n\rightarrow (0,\infty)$, $0< p^{-}\leq p^+\le1<q\leq \infty$.
Fix an integer $d\geq d_{p(\cdot)}\equiv \min\{d\in \mathbb{N}\bigcup\{0\}: p^-(n+d+1)>n\}.$
A function $a$ on $\mathbb{R}^n$ is called a $(p(\cdot),q)$-atom, if there exists a cube $Q$
such that
${\rm supp}\,a\subset Q$;
$\|a\|_{L^q(\mathbb R^n}\leq \frac{|Q|^{1/q}}{\|\chi_{Q}\|_{L^{p(\cdot)}(\mathbb{R}^n)}}$;
$\int_{\mathbb{R}^n} a(x)x^\alpha dx=0\; {\rm for}\;|\alpha| \leq d$.

\begin{theorem}\label{s2t1}\cite{CW,NS,T1}\quad
Let $p(\cdot)\in LH(\mathbb R^n),
0<p^-\le p^+\le1<q<\infty$
and
$$
  \mathcal{A}(\{\lambda_j\}_{j=1}^\infty,\{Q_j\}_{j=1}^\infty)
  =\left\|\left\{\sum_{j}\left(\frac{|\lambda_j|\chi_{Q_j}}{\|\chi_{Q_j}\|_{L^{p(\cdot)}}}\right)^{p^-}
\right\}^{\frac{1}{p^-}}\right\|_{L^{p(\cdot)}}.
$$
If $f\in L^q(\mathbb R^n)\cap H^{p(\cdot)}(\mathbb R^n)$, there is a
sequence of $(p(\cdot),q)-$atoms $\{a_j\}$ and a sequence of scalars $\{\lambda_j\}$
with
\begin{equation*}
  \mathcal{A}(\{\lambda_j\}_{j=1}^\infty,\{Q_j\}_{j=1}^\infty)
  \leq C\|f\|_{{H}^{p(\cdot)}},
\end{equation*}
{\noindent}such that $f=\sum_j\lambda_ja_j$, where the series converges to
$f$ in both $H^{p(\cdot)}(\mathbb R^n)$ and $L^q(\mathbb R^n)$ norms.
Conversely,  if
$$
\mathcal{A}(\{\lambda_j\}_{j=1}^\infty,\{Q_j\}_{j=1}^\infty)<\infty,
$$
then $f=\sum_j\lambda_ja_j$
converges in $\mathcal S_\infty'$, belongs to $H^{p(\cdot)}$ and satisfies
\begin{equation*}
  \|f\|_{{H}^{p(\cdot)}}
  \leq C\mathcal{A}(\{\lambda_j\}_{j=1}^\infty,\{Q_j\}_{j=1}^\infty).
\end{equation*}
\end{theorem}

We also need the following Fefferman-Stein vector valued inequality.
\begin{proposition}\label{s2p1}\cite{CFMP}\quad Let $p(\cdot)\in LH(\mathbb R^n)\cap\mathcal P^0(\mathbb R^n)$.
Then
for any $q>1$, $f=\{f_i\}_{i\in \mathbb{Z}}$, $f_i\in L_{loc}(\mathbb R^n)$, $i\in \mathbb{Z}$
\begin{equation*}
  \|\|\mathbb{M}(f)\|_{l^q}\|_{L^{p(\cdot)}}\leq C\|\|f||_{l^q}\|_{L^{p(\cdot)}},
\end{equation*}
where $\mathbb{M}(f)=\{M(f_i)\}_{i\in\mathbb{Z}}$.
\end{proposition}

Before we prove the theorem, we recall the discrete Calder\'on-type identity.
The well-known discrete Calder\'on identity was first
introduced by Frazier and
Jawerth \cite{FJ1}.  We will
need the following discrete Calder\'on-type identity,
which can be found in \cite{T1} for $\mathbb R^n$ and
was first used in \cite{H,HS} for spaces of homogeneous type.

\begin{lemma}\label{s3l2}\quad Suppose that
$p(\cdot)\in LH(\mathbb R^n)\cap\mathcal P^0(\mathbb R^n)$, $L,N\in \mathbb N$
are large enough. Let $\varphi$ be Schwartz functions with support
on the unit ball satisfying the conditions: for all $\xi\in\mathbb R^n$,
\begin{align*}
\sum_{j\in\mathbb Z}|\widehat\varphi(2^{-j}\xi)|^2=1
\end{align*}
and $\int_{\mathbb R^n}\varphi(x)x^\alpha dx=0$ for
all $0\le|\alpha|\le N$.
Then for all $f\in H^{p(\cdot)}(\mathbb R^n)\cap L^q(\mathbb R^n)$,
$1<q<\infty$, there exists a function
$h\in H^{p(\cdot)}(\mathbb R^n)\cap L^{q}(\mathbb R^n)$ with
\begin{align*}
\|f\|_{L^{q}(\mathbb R^n)}\sim \|h\|_{L^{q}(\mathbb R^n)}
\quad\mbox{and}\quad\|f\|_{H^{p(\cdot)}(\mathbb R^n)}
\sim \|h\|_{H^{p(\cdot)}(\mathbb R^n)}
\end{align*}
such that
\begin{align*}
f(x)&=\sum_{j\in \mathbb Z}\sum_{Q}|Q|\tilde\varphi_j\ast
f(x_Q)\varphi_j(x-x_Q)\\
&=\sum_{j\in \mathbb Z}\sum_{Q}|Q|\varphi_j\ast
h(x_Q)\varphi_j(x-x_Q),
\end{align*}
where $\tilde\varphi\in\mathcal S(\mathbb R^n)$
with $\int_{\mathbb R^n}\tilde\varphi(x)x^\alpha dx=0$, $|\alpha|\le L$ , and
the series converges in both norms of $L^q(\mathbb R^n)$
and $H^{p(\cdot)}(\mathbb R^n)$, and where $x_Q$ is any point in $Q$ and the summation of $Q$ is taken over all dyadic cubes $Q$ with side length $2^{-j-N}$ in $\mathbb{R}^n$ for each $k\in\mathbb Z$ and
a fixed large integer $N$.
\end{lemma}

Then repeating the same argument in
\cite[Proposition 2.3]{T1} and applying Lemma \ref{s3l2} yield that
$$
\|f\|_{H^{p(\cdot)}(\mathbb R^n)}\sim
\left\|\left(\sum_{j\in \mathbb Z}
\sum_{Q}|\varphi_j\ast h(x_Q)|^2\chi_Q\right)^{1/2}
\right\|_{L^{p(\cdot)}(\mathbb R^n)}.
$$

We are ready to prove Theorem \ref{s1th1}.\\
\noindent\textit{\bf Proof of Theorem \ref{s1th1}.}\quad
For $f\in L^2(\mathbb R^n)\cap H^{p(\cdot)}(\mathbb R^n)$
and $\ell\in\mathbb N$, define
$$
\Omega_{\ell}=\{x\in\mathbb R^n:M_{\phi}(f)(x)>2^{\ell}\}
$$
and
$$
B_{\ell}=\left\{Q\;\mbox{is\; a\; dyadic\; cube}:
|Q\cap\Omega_{\ell}|>\frac{1}{2}|Q|,\;|Q\cap\Omega_{\ell+1}|\leq\frac{1}{2}|Q|\right\}.
$$

Let
\begin{align*}
\widetilde\Omega_\ell=\{x\in\mathbb
R^n:M{(\chi_{\Omega_\ell})}(x)>\frac{1}{1000}\}.
\end{align*}
Then $\Omega_\ell\subset\widetilde{\Omega_\ell}$.
By the $L^2$ boundedness of $M$,
$|\widetilde{\Omega_\ell}|\le C|\Omega_\ell|.$

We write $\varphi_Q:=\phi_k,$
if $l(Q)=2^{-k-N}$ and $x_Q$ is any point in $Q$.
Denote
$\tilde Q\in B_\ell$ are maximal dyadic cubes in $B_\ell$, we rewrite
\begin{align*}
\pi_b(f)(x)
&=\sum_k\sum_Q|Q|\varphi_k(x-x_Q)\tilde\varphi_k\ast
b(x_Q)\phi_k\ast f(x_Q)\\
&=\sum_{\ell}\sum_{\tilde Q\in B_\ell}\sum_{Q\subset\tilde Q}|Q|\varphi_Q(x-x_Q)\tilde\varphi_Q\ast
b(x_Q)\phi_Q\ast f(x_Q)\\
&=:\sum_{\ell}\sum_{\tilde Q\in B_\ell}\lambda_{\tilde Q}b_{\tilde Q}(x),
\end{align*}
where
\begin{equation*}
b_{\tilde Q}=\left\{
\begin{aligned}
&0 &, \lambda_{\tilde Q}=0 \\
\frac{1}{\lambda_{\tilde Q}}
\sum_{Q\subset\tilde Q}|Q|\varphi_Q(x-x_Q)\tilde\varphi_Q\ast
&b(x_Q)\phi_Q\ast f(x_Q) &,\lambda_{\tilde Q}\neq 0
\end{aligned}
\right.
\end{equation*}
and
\begin{align*}
\lambda_{\tilde Q}=
C\left\|\left\{\sum_{Q\subset\tilde Q}|\varphi_Q\ast
b(x_Q)|^2|\phi_Q\ast f(x_Q)|^2
\chi_Q\right\}^{\frac{1}{2}}\right\|_{L^{2}(\mathbb R^n)}
|\tilde Q^\ast|^{-1/2}
\|\chi_{\tilde{Q}^\ast}\|_{L^{p(\cdot)}(\mathbb R^n)}.
\end{align*}

By the definition of $b_{\tilde Q}$ and the support of $\varphi$, we have that
$b_{\tilde Q}$ is supported in $\tilde Q^\ast$.
Given $h\in L^{2}(\mathbb R^n)$ with $\|h\|_{L^2(\mathbb R^n)\le 1}$.
By duality and the H\"older's inequalities, we have
\begin{eqnarray*}
&&\left\|\sum_{Q\subset\tilde Q}|Q|\varphi_Q(\cdot-x_Q)\tilde\varphi_Q\ast
b(x_Q)\phi_Q\ast f(x_Q)\right\|_{L^2(\mathbb R^n)}\\
&\le&\sup_{\|h\|_{L^{2}(\mathbb R^n)}\le 1}\bigg<\sum_{Q\subset\tilde Q}|Q|\varphi_Q(\cdot-x_Q)\tilde\varphi_Q\ast
b(x_Q)\phi_Q\ast f(x_Q),h(\cdot)\bigg>\\
&\le&\sup_{\|h\|_{L^{2}(\mathbb R^n)}\le 1}
\left|\sum_{Q\subset\tilde Q}\int(\varphi_Q\ast h(x_Q))(\tilde\varphi_Q\ast
b(x_Q))\phi_Q\ast f(x_Q)\chi_Q(x)dx\right|\\
&\le&\sup_{\|h\|_{L^{2}(\mathbb R^n)}\le 1}
\left\|\left\{\sum_{Q\subset\tilde Q}|\varphi_Q\ast
h(x_Q)|^2\chi_Q\right\}^{\frac{1}{2}}\right\|_{L^{2}(\mathbb R^n)}\\
&&\times\left\|\left\{\sum_{Q\subset\tilde Q}|\varphi_Q\ast
b(x_Q)|^2|\phi_Q\ast f(x_Q)|^2
\chi_Q\right\}^{\frac{1}{2}}\right\|_{L^{2}(\mathbb R^n)}\\
&\le& C
\left\|\left\{\sum_{Q\subset\tilde Q}|
\varphi_Q\ast
b(x_Q)|^2|\phi_Q\ast f(x_Q)|^2
\chi_Q\right\}^{\frac{1}{2}}\right\|_{L^{2}(\mathbb R^n)},
\end{eqnarray*}
where the last inequality follows form the $L^2$ estimates of
the discrete Littlewood-Paley square function estimates.

Then the estimate implies that
\begin{eqnarray*}
\|b_{\tilde Q}\|_{L^2(\mathbb R^n)}&=&\frac{1}{\lambda_{\tilde Q}}
\left\|\sum_{Q\subset\tilde Q}|Q|\varphi_Q(x-x_Q)\tilde\varphi_Q\ast
b(x_Q)\phi_Q\ast f(x_Q)\right\|_{L^2(\mathbb R^n)}\\
&\le &|\tilde Q^\ast|^{1/2}
{\|\chi_{\tilde{Q}^\ast}\|^{-1}_{L^{p(\cdot)}}}.
\end{eqnarray*}

Hence, together with the cancellation conditions of $\varphi$, we
have obtain that $b_{\tilde Q}$ is a $(p(\cdot),2)-$atom.

Therefore, applying the atomic decomposition of $H^{p(\cdot)}(\mathbb R^n)$
in Theorem \ref{s2t1},
\begin{align*}
&\|\pi_b(f)\|_{H^{p(\cdot)}(\mathbb R^n)}\leq
C\mathcal{A}(\{\lambda_\ell\}_{\ell=1}^\infty,\{Q_\ell\}_{\ell=1}^\infty)\\
&\leq C\left\|\left\{\sum_{\ell}\sum_{\tilde Q\in B_\ell}\left(
\frac{|\lambda_{\tilde Q}|\chi_{\tilde Q^\ast}}{\|\chi_{\tilde Q^\ast}\|_{L^{p(\cdot)}}}\right)^{p^-}
\right\}^{\frac{1}{p^-}}\right\|_{L^{p(\cdot)}(\mathbb R^n)}\\
&=C\left\|\left\{\sum_{\ell}\sum_{\tilde Q\in B_\ell}\bigg(
{\bigg\|\bigg\{\sum_{Q\subset\tilde Q}|\varphi_Q\ast
b(x_Q)|^2|\phi_Q\ast f(x_Q)|^2
\chi_Q\bigg\}^{\frac{1}{2}}\bigg\|_{L^{2}(\mathbb R^n)}
|\tilde Q^\ast|^{-1/2}\chi_{\tilde Q^\ast}}\bigg)^{p^-}
\right\}^{\frac{1}{p^-}}\right\|_{L^{p(\cdot)}(\mathbb R^n)}.
\end{align*}

On the other hand,
for any $x\in Q\in B_\ell$, we have
$M\chi_{Q\cap{\tilde\Omega_\ell\setminus\Omega_{\ell+1}}}(x)>\frac{1}{2}$.
Then we have
\begin{align*}
\chi_Q(x)&\le 2M\chi_{Q\cap{\tilde\Omega_\ell\setminus\Omega_{\ell+1}}}(x)\\
&\le 4M^2(\chi_{Q\cap{\tilde\Omega_\ell\setminus\Omega_{\ell+1}}})(x).
\end{align*}

For any $x_Q\in Q\cap\tilde\Omega_\ell\setminus\Omega_{\ell+1}$,
we have $M_\phi(f)(x_Q)\le 2^{\ell+1}$.
By the Fefferman-Stein vector valued inequality,
we have the following estimate:
\begin{align*}
&\left\|\left(
\sum_{Q\subset\tilde Q}|\varphi_Q\ast
b(x_Q)|^2|\phi_Q\ast f(x_Q)|^2
\chi_Q\right\}^{\frac{1}{2}}\right\|^2_{L^2(\mathbb R^n)}\\
&\le\int_{\mathbb R^n}
\sum_{Q\subset\tilde Q}|\varphi_Q\ast
b(x_Q)|^2\big|M_\phi(f)(x_Q)\big|^2
\chi_Q(x)dx\\
&\le C\int_{\mathbb R^n}
\sum_{Q\subset\tilde Q}
|\varphi_Q\ast
b(x_Q)|^2|M_\phi(f)(x_Q)|^2
M^2(\chi_{Q\cap{\tilde\Omega_\ell\setminus\Omega_{\ell+1}}})(x)dx\\
&\le C\int_{\mathbb R^n}
\sum_{Q\subset\tilde Q}
|\varphi_Q\ast
b(x_Q)|^2|M_\phi(f)(x_Q)|^2
\chi_{Q\cap{\tilde\Omega_\ell\setminus\Omega_{\ell+1}}}(x)dx\\
&\le C2^{2\ell}\int_{\tilde Q\cap{\tilde\Omega_\ell\setminus\Omega_{\ell+1}}}
\sum_{Q\subset\tilde Q}\chi_Q(x)
|\varphi_Q\ast b(x_Q)|^2dx\\
&\le C2^{2\ell}
\sum_{Q\subset\tilde Q}|Q|
|\varphi_Q\ast b(x_Q)|^2dx\le C2^{2\ell}|\tilde Q|,
\end{align*}
where the last inequality follows from \cite[page 110, Theorem 4.13]{DH}.

Therefore,
\begin{align*}
&\left\|\left\{\sum_{\ell}\sum_{\tilde Q\in B_\ell}\bigg(
{\bigg\|\bigg\{\sum_{Q\subset\tilde Q}|\varphi_Q\ast
b(x_Q)|^2|\phi_Q\ast f(x_Q)|^2
\chi_Q\bigg\}^{\frac{1}{2}}\bigg\|_{L^{2}(\mathbb R^n)}
|\tilde Q^\ast|^{-1/2}\chi_{\tilde Q^\ast}}\bigg)^{p^-}
\right\}^{\frac{1}{p^-}}\right\|_{L^{p(\cdot)}(\mathbb R^n)}\\
&\le C\left\|\left\{\sum_{\ell}\sum_{\tilde Q\in B_\ell}\left(
2^\ell\chi_{\tilde Q^\ast}\right)^{p^-}
\right\}^{\frac{1}{p^-}}
\right\|_{L^{p(\cdot)}(\mathbb R^n)}.
\end{align*}

Since
\begin{align*}
\chi_{\tilde\Omega_\ell}(x)\le CM^{\frac{2}{p^-}}\chi_{\Omega_{\ell}}(x)
\end{align*}
and
the equivalent characterization of variable Hardy spaces
$$
\|f\|_{H^{p(\cdot)}(\mathbb R^n)}\sim\|M_\phi f\|_{L^{p(\cdot)}(\mathbb R^n)},
$$
then by using
the Fefferman-Stein vector valued inequality in Proposition \ref{s2p1}, we yields that
\begin{align*}
&\left\|\left\{\sum_{\ell}\sum_{\tilde Q\in B_\ell}\left(
2^\ell\chi_{\tilde Q^\ast}\right)^{p^-}
\right\}^{\frac{1}{p^-}}
\right\|_{L^{p(\cdot)}(\mathbb R^n)}\\
&\le C\left\|\left\{\sum_{\ell}\left(
2^\ell M^{\frac{2}{p^-}}\chi_{\Omega_{\ell}}\right)^{p^-}
\right\}^{\frac{1}{p^-}}
\right\|_{L^{p(\cdot)}(\mathbb R^n)}
= C\left\|\left\{\sum_{\ell}
2^{\ell p^-}M^{2}\chi_{\Omega_{\ell}}
\right\}^{\frac{1}{2}}
\right\|^{\frac{2}{p^-}}_{L^{\frac{2p(\cdot)}{p^-}}(\mathbb R^n)}\\
&\le C\left\|\left\{\sum_{\ell}
2^{\ell p^-}\chi_{\Omega_{\ell}}^2
\right\}^{\frac{1}{2}}
\right\|^{\frac{2}{p^-}}_{L^{\frac{2p(\cdot)}{p^-}}(\mathbb R^n)}
\le C\left\|\left\{\sum_{\ell}\left(
2^\ell\chi_{ \Omega_\ell}\right)^{p^-}
\right\}^{\frac{1}{p^-}}
\right\|_{L^{p(\cdot)}(\mathbb R^n)}\\
&\le C\left\|\left\{\sum_{\ell}\left(
2^\ell\chi_{ \Omega_\ell\setminus\Omega_{\ell+1}}\right)^{p^-}
\right\}^{\frac{1}{p^-}}
\right\|_{L^{p(\cdot)}(\mathbb R^n)}\\
&= C\inf\left\{\lambda>0:\sum_\ell\int_{{ \Omega_\ell\setminus\Omega_{\ell+1}}}
\left(\frac{2^\ell}{\lambda}
\right)^{p(x)}dx\le 1\right\}\\
&\le C\inf\left\{\lambda>0:\int_{\mathbb R^n}
\left(\frac{M_\phi f(x)}{\lambda}
\right)^{p(x)}dx\le 1\right\}\leq C\|f\|_{{H}^{p(\cdot)}(\mathbb R^n)}.
\end{align*}

If we combine these estimates we get that
$$
\|\pi_b(f)\|_{H^{p(\cdot)}(\mathbb R^n)}
\le C\|f\|_{H^{p(\cdot)}(\mathbb R^n)}
$$
for $f\in L^2(\mathbb R^n)\cap H^{p(\cdot)}(\mathbb R^n)$.
Since that
$L^2(\mathbb R^n)\cap H^{p(\cdot)}(\mathbb R^n)$
is dense in
${H}^{p(\cdot)}(\mathbb R^n)$,
then
by the density argument
$\pi_b$ can be extended to a bounded operator
on $H^{p(\cdot)}(\mathbb R^n)$.

Now we prove that $\pi_b$ is bounded on $L^2(\mathbb R^n)$.
By applying H\"older's inequality and Carleson's condition,
\begin{eqnarray*}
&&\|\pi_b(f)\|_{L^2(\mathbb R^n)}
=\sup_{\|g\|_{L^2(\mathbb R^n)}\le 1}|\big<\pi_b(f),g\big>|\\
&=&\sup_{\|g\|_{L^{2}(\mathbb R^n)}\le 1}\bigg<\sum_k\sum_Q|Q|\varphi_Q(\cdot-x_Q)\tilde\varphi_Q\ast
b(x_Q)\phi_Q\ast f(x_Q),g(\cdot)\bigg>\\
&\le&\sup_{\|g\|_{L^{2}(\mathbb R^n)}\le 1}
\left|\sum_k\sum_Q|Q|(\varphi_Q\ast g(x_Q))(\tilde\varphi_Q\ast
b(x_Q))\phi_Q\ast f(x_Q)\chi_Q(x)\right|\\
&\le&\sup_{\|h\|_{L^{2}(\mathbb R^n)}\le 1}
\left\{\sum_{k}\sum_{Q}
|Q||\varphi_Q\ast
g(x_Q)|^2\chi_Q\right\}^{\frac{1}{2}}\\
&&\times\left(\sum_{\ell}\sum_{\tilde Q\in B_\ell}\sum_{Q\subset\tilde Q}
|Q||\varphi_Q\ast
b(x_Q)|^2|\phi_Q\ast f(x_Q)|^2
\chi_Q\right)^{\frac{1}{2}}\\
&\le& C
\left(\sum_{\ell}\sum_{\tilde Q\in B_\ell}\sum_{Q\subset\tilde Q}
|Q||\varphi_Q\ast
b(x_Q)|^2|\phi_Q\ast f(x_Q)|^2
\right)^{\frac{1}{2}}\\
&\le& C
\left(\sum_{\ell}2^{2(\ell+1)}\sum_{\tilde Q\in B_\ell}\sum_{Q\subset\tilde Q}
|Q||\varphi_Q\ast
b(x_Q)|^2
\right)^{\frac{1}{2}}\\
&\le& C
\bigg(\sum_{\ell}2^{2(\ell+1)}|\tilde\Omega_\ell|
\bigg)^{\frac{1}{2}}\le C\|f\|_{L^2(\mathbb R^n)},
\end{eqnarray*}
for any $f\in L^2(\mathbb R^n)$.
Thus, $\pi_b$ is bounded on $L^2(\mathbb R^n)$
and $H^{p(\cdot)}(\mathbb R^n).$
Note that $\sup\limits_{k}|\phi_k\ast f(x)|\le CM(f)(x)$.
Then for $q>1$, $\sup\limits_{k}|\phi_k\ast f(x)|\in L^q(\mathbb R^n)$.
Since
$$
\lim_{t\rightarrow0}\|\phi_t\ast f-f\|_{L^q(\mathbb R^n)}=0,
$$
there is a sequence of $t_j\rightarrow 0$ such that
$\lim_{t_j\rightarrow0}\phi_{t_j}\ast f(x)=f(x)$
for a.e. $x\in\mathbb R^n$. Then for $p(\cdot)\in \mathcal P^0$,
$$
\|f\|_{L^{p(\cdot)}(\mathbb R^n)}\le\lim_{t_j\rightarrow0}
\|\phi_{t_j}\ast f\|_{L^{p(\cdot)}(\mathbb R^n)}
$$
and thus
$$
\|f\|_{L^{p(\cdot)}(\mathbb R^n)}
\le C\|f\|_{H^{p(\cdot)}(\mathbb R^n)}.
$$

Therefore, we get that
$$
\|\pi_b(f)\|_{L^{p(\cdot)}(\mathbb R^n)}
\le C\|\pi_b(f)\|_{H^{p(\cdot)}(\mathbb R^n)}
\le C\|f\|_{H^{p(\cdot)}(\mathbb R^n)}
$$
for $f\in L^2(\mathbb R^n)\cap H^{p(\cdot)}(\mathbb R^n)$.

Similarly,
by the density argument
$\pi_b$ can be extended to a bounded operator
from $H^{p(\cdot)}(\mathbb R^n)$ to $H^{p(\cdot)}(\mathbb R^n)$.
Thus we complete the proof of Theorem \ref{s1th1}.
$\hfill\Box$

\section{Proof of Theorem \ref{s1th2} and Theorem \ref{s1th3}}

In this section, we will prove the $H^{p(\cdot)}$ boundedness of
Calder\'on-Zygmund operators. For the proof we first need the
discrete Littlewood-Paley characterizations for $H^{p(\cdot)}$
in \cite[Proposition 2.3]{T1}.
Let $\psi\in\mathcal S(\mathbb R^n)$
satisfy
\begin{align*}
\mbox{supp} \widehat\psi\subset\{\xi\in\mathbb R^n:1/2\le|\xi|\le2\}
\end{align*}
and
\begin{align*}
\sum_{j\in\mathbb Z}|\widehat\psi(2^{-j}\xi)|^2=1\quad \mbox{for\;all}\;\xi\in\mathbb{R}^n\setminus\{0\}.
\end{align*}
Denote by $\mathcal S_\infty(\mathbb R^n)$ the functions
$f\in\mathcal S(\mathbb R^n)$ satisfying
$\int_{\mathbb R^n}f(x)x^\alpha dx=0$ for all muti-indices $\alpha\in \mathbb Z_+^n:=(\{0,1,2,\cdots\})^n$
and $\mathcal S'_\infty(\mathbb R^n)$ its topological dual space. For $f\in \mathcal S'_\infty(\mathbb R^n)$,
we recall the definition of the Littlewood-Paley-Stein square function
\begin{align*}
\mathcal{G}(f)(x):=\left(\sum_{j\in \mathbb Z}|\psi_j\ast f(x)|^2\right)^{1/2},
\end{align*}
and the discrete Littlewood-Paley-Stein square function
\begin{align*}
\mathcal{G}^d(f)(x):=\left(\sum_{j\in \mathbb Z}
\sum_{\mathbf k\in\mathbb Z^n}|\psi_j\ast f(2^{-j}\mathbf k)|^2\chi_Q(x)\right)^{1/2},
\end{align*}
where $Q$ denote dyadic cubes in $\mathbb R^n$ with side-lengths $2^{-j}$ and the
lower left-corners of $Q$ are $2^{-j}\mathbf k$.
If $p(\cdot)\in LH(\mathbb R^n)\cap \mathcal P^0(\mathbb R^n)$,
$$
\|f\|_{H^{p(\cdot)}(\mathbb R^n)}\sim
\|\mathcal G(f)\|_{L^{p(\cdot)}(\mathbb R^n)}\sim \|\mathcal G^d(f)\|_{L^{p(\cdot)}(\mathbb R^n)}.
$$
Then by discrete Calder\'on-type identity, the
discrete Littlewood-Paley characterizations for $H^{p(\cdot)}$
and almost orthogonal estimates,
we show the following theorem.

\begin{theorem}\label{s3th1}\quad
Suppose that $p(\cdot)\in LH(\mathbb R^n)$ and $T\in CZO(\epsilon),$ where $0<\epsilon\le 1$. If $T1=T^\ast1=0$,
then $T$ extends to a continuous operator on $H^{p(\cdot)}(\mathbb R^n)$,$\frac{n}{n+\epsilon}<p^-\le p^+\le1$.
\end{theorem}

Now we are ready to prove Theorem \ref{s3th1}.\\

\noindent\textit{\bf Proof of Theorem \ref{s3th1}}\quad Given $q>1$,
since that the subspace $H^{p(\cdot)}(\mathbb R^n)\cap L^q(\mathbb R^n)$
is dense in $H^{p(\cdot)}(\mathbb R^n)$,
we only need to prove that $T$ is bounded from
$H^{p(\cdot)}(\mathbb R^n)\cap L^q(\mathbb R^n)$
to $H^{p(\cdot)}(\mathbb R^n)$.
For any $f\in H^{p(\cdot)}(\mathbb R^n)\cap L^q(\mathbb R^n)$
and $T\in CZO(\epsilon)$, we have $Tf\in L^q(\mathbb R^n)$.
Then by Lemma \ref{s3l2},
\begin{align*}
f(x)=\sum_{j'\in \mathbb Z}\sum_{Q'}|Q'|\varphi_{j'}\ast
h(x_{Q'})\varphi_{j'}(x-x_{Q'})
\end{align*}
and
\begin{align*}
\varphi_j\ast Tf(x_Q)=\sum_{j'\in \mathbb Z}\sum_{Q'}|Q'|\varphi_{j'}\ast
h(x_{Q'})\varphi_jT\varphi_{j'}(x_Q,x_{Q'}),
\end{align*}
where
\begin{align*}
\begin{split}
\varphi_jT\varphi_{j'}(x_Q,x_{Q'})
=\int\int\varphi_{j}(x_Q-u)K(u,v)\varphi_{j'}(v-x_{Q'})dudv
\end{split}
\end{align*}
and
$h\in H^{p(\cdot)}(\mathbb R^n)\cap L^{q}(\mathbb R^n)$ with
\begin{align*}
\|f\|_{L^{q}(\mathbb R^n)}\sim \|h\|_{L^{q}(\mathbb R^n)}
\quad\mbox{and}\quad\|f\|_{H^{p(\cdot)}(\mathbb R^n)}
\sim \|h\|_{H^{p(\cdot)}(\mathbb R^n)}.
\end{align*}

We claim that

\begin{align}\label{s3i1}
\begin{split}
|\varphi_jT\varphi_{j'}(x_Q,x_{Q'})|
&\leq C2^{-|j-j'|\epsilon}
\frac{2^{-(j\wedge j')\epsilon}}{[2^{-(j\wedge j')}+|x_Q-x_{Q'}|]^{n+\epsilon}}.
\end{split}
\end{align}

To prove (\ref{s3i1}), we consider the four cases:
(1) $j> j'$ and $|x_Q-x_{Q'}|\leq 5\;2^{-j'}$;
(2) $j> j'$ and $|x_Q-x_{Q'}|\ge 5\;2^{-j'}$;
(3) $j\leq j'$ and $|x_Q-x_{Q'}|\leq 5\;2^{-j'}$;
(4) $j\leq j'$ and $|x_Q-x_{Q'}|\ge 5\;2^{-j'}$.
The idea we used here
comes from \cite{HS1}.\\
In Case (1),
since $T1=0$, we have
\begin{align*}
\begin{split}
\varphi_jT\varphi_{j'}(x_Q,x_{Q'})=&
\int\int\varphi_{j}(x_Q-u)K(u,v)\varphi_{j'}(v-x_{Q'})dudv\\
=&\int\int\varphi_{j}(x_Q-u)K(u,v)(\varphi_{j'}(v-x_{Q'})-\varphi_{j'}(x_Q-x_{Q'}))dudv.
\end{split}
\end{align*}

Choose a smooth function $\eta_0$ such that $\mbox{supp}\;\eta_0\subset
\{x:|x|\leq 6\}$ and let $\eta_0=1$ when $|x|\leq 2$. Set $\eta_1=
1-\eta_0$. Then we get that
\begin{align*}
\begin{split}
&|\varphi_jT\varphi_{j'}(x_Q,x_{Q'})|\\
=&\bigg|\int\int\varphi_j(x_Q-u)K(u,v)(\varphi_{j'}(v-x_{Q'})-\varphi_{j'}(x_Q-x_{Q'}))
\eta_0(2^{j}(v-x_Q))dudv\bigg|\\
&+\bigg|\int\int\varphi_j(x_Q-u)K(u,v)(\varphi_{j'}(v-x_{Q'})-\varphi_{j'}(x_Q-x_{Q'}))
\eta_1(2^{j}(v-x_Q))dudv\bigg|\\
=&I+II.
\end{split}
\end{align*}

For $I$, we denote
$\tilde\varphi(v)=(\varphi_{j'}(v-x_{Q'})-\varphi_{j'}(x_Q-x_{Q'}))
\eta_0(2^{j}(v-x_Q))$
and $\omega(u)=\varphi_j(x_Q-u)$. Applying H\"older's inequality and the
$L^2(\mathbb R^n)$ boundedness of
$T$ yield that
\begin{align*}
\begin{split}
&I= |\big<T\tilde\varphi,\omega\big>|\leq \|T\tilde\varphi\|_{L^2(\mathbb R^n)}
\|\omega\|_{L^2(\mathbb R^n)}\\
&\le C2^{j'-j}2^{-j'n}.
\end{split}
\end{align*}

We now deal with the term $II$.
By the cancellation condition of $\varphi$, we get that
\begin{align*}
\begin{split}
&II
=\bigg|\int\int\varphi_j(x_Q-u)[K(u,v)-K(x_Q,v)]
\\&\times(\varphi_{j'}(v-x_{Q'})-\varphi_{j'}(x_Q-x_{Q'}))
\eta_1(2^{j}(v-x_Q))dudv\bigg|\\
&\le C2^{-(j-j')\epsilon}2^{-j'n}.
\end{split}
\end{align*}

In Case (2), observe that $|x_Q-x_{Q'}|\sim|u-v|$.
The smoothness condition on the kernel $K(u,v)$ implies
\begin{align*}
\begin{split}
|\varphi_jT\varphi_{j'}(x_Q,x_{Q'})|\leq&
\bigg|\int\int\varphi_{j}(x_Q-u)K(u,v)\varphi_{j'}(v-x_{Q'})dudv\bigg|\\
=&\int\int|\varphi_{j}(x_Q-u)||K(u,v)-K(x_Q,v)||\varphi_{j'}(v-x_{Q'})|dudv\\
\le& C\frac{2^{-j\epsilon}}{|x_Q-x_{Q'}|^{n+\epsilon}}.
\end{split}
\end{align*}

The other cases are similar to Case (1) and Case (2).
So we prove the claim.

Then using the claim and the Fefferman-Stein vector-valued maximal
inequality, we have
\begin{align*}
\begin{split}
&\|Tf\|_{H^{p(\cdot)}(\mathbb R^n)}
\sim
\left\|\left(\sum_{j\in \mathbb Z}
\sum_{Q}|\varphi_j\ast Tf(x_Q)|^2\chi_Q\right)^{1/2}
\right\|_{L^{p(\cdot)}(\mathbb R^n)}\\
&\leq
C\left\|\left(\sum_{j\in \mathbb Z}
\sum_{Q}\big|\sum_{j'\in \mathbb Z}\sum_{Q'}|Q'|\varphi_{j'}\ast
h(x_{Q'})\varphi_jT\varphi_{j'}(x_Q,x_{Q'})\big|^2\chi_Q\right)^{1/2}
\right\|_{L^{p(\cdot)}(\mathbb R^n)}\\
&\leq
C\left\|\Bigg(\sum_{j'\in \mathbb Z}
   \left\{M\left(\sum_{Q}
|(\varphi_{j'}\ast h)x_{Q'}|^2\chi_{Q'}
\right)^{\frac{\delta}{2}}\right\}
^{\frac{2}{\delta}}\Bigg)^{\frac{\delta}{2}}
\right\|^\frac{1}{\delta}_{L^{\frac{p(\cdot)}{\delta}}}\\
&\leq C\left\|\left(\sum_{j'\in \mathbb Z}
\sum_{Q'}|\varphi_{j'}\ast h(x_{Q'})|^2\chi_{Q'}\right)^{1/2}
\right\|_{L^{p(\cdot)}(\mathbb R^n)}\sim
\|f\|_{H^{p(\cdot)}(\mathbb R^n)},
\end{split}
\end{align*}
where the second inequality follows from the lemma in \cite[pages 147-148]{FJ2} and
$\frac{n}{n+\epsilon}<\delta<p^-$.

Thus, we have completed the proof of Theorem \ref{s3th1}.
$\hfill\Box$

We now turn to the

\noindent\textit{\bf Proof of Theorem \ref{s1th2}:}\quad
Let $p(\cdot)\in LH(\mathbb R^n)$ and $T\in CZO(\epsilon),$ where $0<\epsilon\le 1$.
First, by using para-product operator we prove that $T$ is bounded on $H^{p(\cdot)}(\mathbb R^n)$,$\frac{n}{n+\epsilon}<p^-\le p^+\le1$
if $T^\ast1=0$. From Lemma \ref{s3l2}, we obtain
\begin{align*}
\pi_b(1)(x)&=\sum_k\sum_Q|Q|\varphi_k(x-x_Q)\tilde\varphi_k\ast
b(x_Q)\phi_k\ast 1(x_Q)\\
&=\sum_{k}\sum_{Q}|Q|\varphi_k(x-x_Q)\tilde\varphi_k\ast
b(x_Q)=b(x)
\end{align*}
and
\begin{align*}
\pi^\ast_b(1)(x)&=\sum_k\sum_Q|Q|\varphi_k\ast 1(x_Q)\tilde\varphi_k\ast
b(x_Q)\phi_k(x_Q-x)\\
&=\sum_k\sum_Q|Q|\bigg(\int_{\mathbb R^n}\varphi_k(x)dx\bigg)
\tilde\varphi_k\ast
b(x_Q)\phi_k(x_Q-x)=0.
\end{align*}

In the proof of Theorem \ref{s1th1}, we have showed that
$\pi_b$ is bounded on $L^2(\mathbb R^n)$.
Note that the kernel of $\pi_b$ is the function
$$
K_b(x,y)=\sum_k\sum_Q|Q|\varphi_k(x-x_Q)\tilde\varphi_k\ast
b(x_Q)\phi_k(x_Q-y)
$$
fulfilling the all conditions (\ref{s1i1}), (\ref{s1i2}) and (\ref{s1i3}).
In fact, applying the size condition of $\phi,\varphi$ and
$\tilde\varphi\in H^1(\mathbb R^n)$ yield
\begin{align*}
|K_b(x,y)|&=
\bigg|\sum_k\sum_Q|Q|\varphi_k(x-x_Q)\tilde\varphi_k\ast
b(x_Q)\phi_k(x_Q-y)\bigg|\\
&\leq
C\sum_k\sum_Q|Q||\varphi_k(x-x_Q)|\|\tilde\varphi_k\|_{H^1(\mathbb R^n)}
\|b\|_{BMO(\mathbb R^n)}\phi_k(x_Q-y)|\\
&\leq C\|b\|_{BMO}
\sum_k\sum_Q|Q||\varphi_k(x-x_Q)|
\frac{2^{-k}}{(2^{-k}+|x_Q-y|)^{n+1}}\\
&\leq C
\sum_k\
\frac{2^{-k}}{(2^{-k}+|x-y|)^{n+1}}
\leq C
\frac{1}{|x-y|^{n}}.
\end{align*}
Repeating the similar argument, the kernel
$K_b(x,y)$ also
satisfies the conditions (\ref{s1i2}) and (\ref{s1i3})
and then $\pi_b$ is a Calder\'on-Zygmund operator.
We define the new operator $\tilde T$ by
$$
\tilde T=T-\pi_{T1}.
$$
By Theorem \ref{s3th1}, $\tilde T$ is bounded on
$H^{p(\cdot)}(\mathbb R^n)$ since
$\tilde T1=T1-\pi_{T1}1=0$ and
$\tilde T^\ast1=T^\ast1-\pi^\ast_{T1}1=0$.
By Theorem \ref{s1th1}, $\pi_b$ is bounded on
$H^{p(\cdot)}(\mathbb R^n)$ and then $T$ is bounded on
$H^{p(\cdot)}(\mathbb R^n)$.  Therefore, the condition
$T^\ast1=0$ is sufficient. On the other hand, this condition
is obviously necessary. Indeed,
for any $\varphi\in\mathcal D_0$, then
$\varphi\in H^{1}(\mathbb R^n)$ and $T\varphi\in H^{1}(\mathbb R^n)$. Thus,
$\int T\varphi(x)dx=0$.
Therefore, we conclude the proof of Theorem \ref{s1th2}.
$\hfill\Box$

\noindent\textit{\bf Proof of Theorem \ref{s1th3}}\quad
The proof is similar to the ones of Theorem \ref{s3th1}
and Theorem \ref{s1th2}. We only need to
observe that $\varphi$ has the zero vanishing moment up to order $N$
and $T^\ast(x^\beta)=0$ with $|\beta|\le N$.
Then by repeating the similar argument together with
using Taylor expansion of the kernel of $T$ and
the high order moment condition of $\varphi$,
then we can establish the desired
almost orthogonal estimate.
On the other hand, for any bounded, compactly supported function
$\varphi$ fulfilling the moment condition $\int x^\alpha \varphi(x)dx=0$
for all $|\alpha|\le N$, then we also have
$\int x^\alpha T\varphi(x)dx=0$. The integral is well defined since
$T\varphi(x)=\mathcal O(|x|^{-n-N})$ at infinity.
This means
$T^\ast x^\alpha=0.$
$\hfill\Box$

{\bf Acknowledgments.}
The project is sponsored by Natural Science Foundation of Jiangsu Province of China (grant no. BK20180734), Natural Science Research of Jiangsu Higher Education Institutions of China (grant no. 18KJB110022) and Nanjing University of Posts and Telecommunications Science Foundation (grant no. NY217151).

\bibliographystyle{amsplain}

\end{document}